\documentclass[11pt]{article}
\usepackage{a4wide}

\newcommand{\ba}{\begin{array}}
\newcommand{\ea}{\end{array}}
\newcommand{\be}{\begin{equation}}
\newcommand{\ee}{\end{equation}}
\newcommand{\la}{\label}
\newcommand{\bea}{\begin{eqnarray}}
\newcommand{\eea}{\end{eqnarray}}
\newcommand{\ch}{\choose}
\renewcommand{\a}{\alpha}
\renewcommand{\b}{\beta}
\renewcommand{\c}{\gamma}
\newcommand{\G}{\Gamma}
\renewcommand{\L}{L_n^{(\a)}(x)}
\newcommand{\KL}{L_n^{\a,M}(x)}
\newcommand{\gL}{L_n^{\a,M,N}(x)}
\renewcommand{\P}{P_n^{(\a,\b)}(x)}
\newcommand{\GP}{P_n^{\a,\b,M,N}(x)}
\newcommand{\SP}{P_n^{(\a,\a)}(x)}
\newcommand{\SGP}{P_n^{\a,\a,M,M}(x)}
\renewcommand{\l}{\left}
\renewcommand{\r}{\right}
\newcommand{\set}[1]{\left\{#1\right\}_{n=0}^{\infty}}
\newcommand{\hyp}[5]{\mbox{}_{#1}F_{#2}
\left(\left.\begin{array}{c}#3\\#4\end{array}\right|#5\right)}
\newcommand{\ndots}{n=0,1,2,\ldots}
\newcommand{\n}{\nonumber}
\newcommand{\nn}{\nonumber \\}
\newcommand{\ds}{\displaystyle}
\newcounter{stelling}
\newcommand{\st}[1]{\par\vspace{0.5cm}\refstepcounter{stelling}
{\bf Theorem \thestelling.} {\sl #1}\par\vspace{0.5cm}}

\begin{document}

\begin{center}
{\Large THE SEARCH FOR\\ DIFFERENTIAL EQUATIONS\\ FOR CERTAIN SETS OF\\
\vspace{2mm} ORTHOGONAL POLYNOMIALS}

\vspace{1cm}

{\large Roelof Koekoek}

\vspace{1cm}

Delft University of Technology\\ Faculty of Technical Mathematics and Informatics\\
Mekelweg 4\\ 2628 CD Delft\\ The Netherlands\\ e-mail : koekoek@twi.tudelft.nl
\end{center}

\vspace{2cm}

\begin{abstract}
We look for differential equations of the form
$$\sum_{i=0}^{\infty}c_i(x)y^{(i)}(x)=0,$$
where the coefficients $\l\{c_i(x)\r\}_{i=0}^{\infty}$ are continuous
functions on the real line and where $\l\{c_i(x)\r\}_{i=1}^{\infty}$
are independent of $n$, for the generalized Jacobi polynomials $\set{\GP}$
and for generalized Laguerre polynomials $\set{\gL}$ which are
orthogonal with respect to an inner product of Sobolev type.

We use a method involving computeralgebra packages like Maple
and Mathematica and we will give some preliminary results.
\end{abstract}

\vspace{1cm}


\newpage

\section{Introduction}

In this paper we consider the polynomials $\set{\GP}$ (see \cite{Koorn})
which are orthogonal on the interval $[-1,1]$ with respect to the weight
function
$$\frac{\G(\a+\b+2)}{2^{\a+\b+1}\G(\a+1)\G(\b+1)}(1-x)^{\a}(1+x)^{\b}
+M\delta(x+1)+N\delta(x-1),$$
where $\a>-1$, $\b>-1$, $M\ge 0$ and $N\ge 0$.

We also consider the polynomials $\set{\gL}$ (see \cite{SIAM}) which are
orthogonal with respect to the Sobolev inner product
$$<f,g>\;=\frac{1}{\G(\a+1)}\int\limits_0^{\infty}x^{\a}e^{-x}f(x)g(x)dx+
Mf(0)g(0)+Nf'(0)g'(0),$$
where $\a>-1$, $M\ge 0$ and $N\ge 0$.

We are looking for differential equations of the form
$$\sum_{i=0}^{\infty}c_i(x)y^{(i)}(x)=0,$$
satisfied by these sets of orthogonal polynomials, where the coefficients
$\l\{c_i(x)\r\}_{i=0}^{\infty}$ are continuous
functions on the real line and $\l\{c_i(x)\r\}_{i=1}^{\infty}$ are
independent of the degree $n$.

In \cite{DV} J.~Koekoek and R.~Koekoek found a differential equation of the
above type for the generalized Laguerre polynomials $\set{\KL}$ (see \cite{Koorn})
for all $\a>-1$. These polynomials form a special case of the above mentioned
Sobolev Laguerre polynomials, since $L_n^{\a,M,0}(x)=\KL$.
For more details the reader is referred to \cite{DV} and section 2 of this paper
for more results concerning the coefficients of this differential equation.

The differential equation found in \cite{DV} was computed by hand, without the
help of computers. This is nearly impossible for these other cases. We need
computers to handle the very huge expressions we have to deal with.

The method we used to find the results listed in this paper is explained in
\cite{Rap}. We refer to this report for more details.

\section{The infinite order Laguerre differential equation}

In \cite{DV} the following theorem was proved without the use of computers~:

\st{For $M>0$ the polynomials $\set{\KL}$ satisfy a unique differential equation
of the form
\be\la{dv}M\sum_{i=0}^{\infty}a_i(x)y^{(i)}(x)+xy''(x)+(\a+1-x)y'(x)+ny(x)=0,\ee
where $\l\{a_i(x)\r\}_{i=0}^{\infty}$ are continuous functions on the real
line and $\l\{a_i(x)\r\}_{i=1}^{\infty}$ are independent of $n$.

Moreover, the functions $\l\{a_i(x)\r\}_{i=0}^{\infty}$ are polynomials given by
\be\la{coeff}\l\{\ba{l}\ds a_0(x)={n+\a+1 \ch n-1}\\ \\
\ds a_i(x)=\frac{1}{i!}\sum_{j=1}^i(-1)^{i+j+1}{\a+1 \ch j-1}{\a+2 \ch i-j}(\a+3)_{i-j}x^j,
\;i=1,2,3,\ldots.\ea\r.\ee}

Later we discovered that the coefficients $\l\{a_i(x)\r\}_{i=1}^{\infty}$
have the following interesting property.

\st{The coefficients $\l\{a_i(x)\r\}_{i=1}^{\infty}$ of the differential
equation given by (\ref{dv}) and (\ref{coeff}) satisfy
$$\sum_{i=1}^{\infty}a_i(x)=-\frac{\sin\pi\a}{\pi}\frac{x}{(\a+2)(\a+3)}
\hyp{1}{1}{1}{\a+4}{-x},\;\a>-1.$$

For nonnegative integer values of $\a$ we have~:
$$\sum_{i=k}^{\infty}{i \ch k}a_i(x)=(-1)^{\a+k}a_k(-x),\;k=1,2,3,\ldots.$$}

Note that this theorem implies for nonnegative integer values of $\a$~:
$$\sum_{i=1}^{\infty}a_i(x)=0\;\mbox{ and }\;\sum_{i=1}^{\infty}ia_i(x)=(-1)^{\a+1}x.$$

The proof of this theorem can be found in \cite{Rap}.

\section{Some preliminary results for the Sobolev Laguerre polynomials}

In this section we look for a differential equation of the form
\bea\la{Sobdv} & &M\sum_{i=0}^{\infty}a_i(x)y^{(i)}(x)+
N\sum_{i=0}^{\infty}b_i(x)y^{(i)}(x)+{}\nn
& &\hspace{1cm}{}+MN\sum_{i=0}^{\infty}c_i(x)y^{(i)}(x)+
xy''(x)+(\a+1-x)y'(x)+ny(x)=0\eea
for the polynomials
$$y(x)=\gL=A_0\L+A_1\frac{d}{dx}\L+A_2\frac{d^2}{dx^2}\L,$$
where the coefficients $A_0$, $A_1$ and $A_2$ are defined by
\be\la{AAA}\l\{\ba{l}\ds A_0=1+M{n+\a \ch n-1}+\frac{n(\a+2)-(\a+1)}{(\a+1)(\a+3)}
N{n+\a \ch n-2}+{}\\
\ds\hspace{5cm}{}+\frac{MN}{(\a+1)(\a+2)}{n+\a \ch n-1}{n+\a+1 \ch n-2}\\  \\
\ds A_1=M{n+\a \ch n}+\frac{(n-1)}{(\a+1)}N{n+\a \ch n-1}+
\frac{2MN}{(\a+1)^2}{n+\a \ch n}{n+\a+1 \ch n-2}\\  \\
\ds A_2=\frac{N}{(\a+1)}{n+\a \ch n-1}+\frac{MN}{(\a+1)^2}{n+\a \ch n}
{n+\a+1 \ch n-1}.\ea\r.\ee
For details concerning these generalized Laguerre polynomials and their
definition the reader is referred to \cite{SIAM} and \cite{Thesis}.

Of course, since $L_n^{\a,M,0}(x)=\KL$, the coefficients $\l\{a_i(x)\r\}_{i=0}^{\infty}$
are given by (\ref{coeff}).

Although the general form is still an open problem so far, we know that the
differential equation given by (\ref{Sobdv}) is not unique as in the case of
the differential equation (\ref{dv}).
This is explained by the following theorem~:

\st{The polynomials $\l\{\gL\r\}_{n=1}^{\infty}$ satisfy the following
infinite order differential equation~:
$$\sum_{i=0}^{\infty}b_i^*(\a,x)y^{(i)}(x)+
M\sum_{i=0}^{\infty}c_i^*(\a,x)y^{(i)}(x)=0,$$
where
$$b_i^*(\a,x)=\frac{1}{i!}\sum_{j=0}^i(-1)^j{i \ch j}
(\a+1)_{i-j}x^j,\;i=0,1,2,\ldots$$
and
$$c_i^*(\a,x)=\frac{(-1)^i}{i!}x^i,\;i=0,1,2,\ldots.$$}

The proof is very easy and is based on the observation that
$$\sum_{i=0}^{\infty}b_i^*(\a,x)D^{i+k}\L
=\frac{(-n)_k}{n\G(k)},\;n\ge 1,\;k=0,1,2$$
and
$$\sum_{i=0}^{\infty}c_i^*(\a,x)D^{i+k}\L
={n+\a \ch n}\frac{(-n)_k}{(\a+1)_k},\;k=0,1,2.$$

By using these formulas and the definition (\ref{AAA}) of the coefficients
$A_0$, $A_1$ and $A_2$ the result follows immediately.

Again the reader is referred to the report \cite{Rap} for more details.

\vspace{5mm}

In the special cases $\a=0$, $\a=1$ and $\a=2$ differential equations of the
form (\ref{Sobdv}) are found explicitly. In these
three special cases of integer values of the parameter $\a$ we find a linear
differential equation of formal order $4\a+10$. By formal order we mean that
for special cases ($M=0$ or $N=0$) the true order might be lower.
We will not give all results here, but we again refer to \cite{Rap} for
all details.

As an example we give the result in the case that $\a=0$.
In that case we have found the following linear differential equation of
formal order $10$~:
\bea & &\frac{1}{60}MNx^5y^{(10)}(x)+\frac{1}{12}MN(5x^4-x^5)y^{(9)}(x)+{}\nn
& &{}+\l[\frac{1}{24}MN(72x^3-45x^4+4x^5)-\frac{1}{12}Nx^4\r]y^{(8)}(x)+{}\nn
& &{}+\l[\frac{1}{6}MN(36x^2-72x^3+20x^4-x^5)-\frac{1}{3}N(4x^3-x^4)\r]
y^{(7)}(x)+{}\nn
& &{}+\l[\frac{1}{12}MN(-252x^2+224x^3-35x^4+x^5)
+\frac{1}{2}N(-10x^2+9x^3-x^4)\r]y^{(6)}(x)+{}\nn
& &{}+\l[\frac{1}{60}MN(1680x^2-840x^3+75x^4-x^5)
+\frac{1}{6}N(-12x+81x^2-33x^3+2x^4)\r]y^{(5)}(x)+{}\nn
& &{}+\l[\frac{1}{24}MN(-420x^2+120x^3-5x^4)+{}\r.\nn
& &\hspace{3cm}\l.{}+\frac{1}{12}N(24+36x-150x^2+34x^3-x^4)
-\frac{1}{2}Mx^2\r]y^{(4)}(x)+{}\nn
& &{}+\l[\frac{1}{3}MN(15x^2-2x^3)
+\frac{1}{2}N(-6-2x+9x^2-x^3)+M(-2x+x^2)\r]y^{(3)}(x)+{}\nn
& &{}+\l[-\frac{1}{2}MNx^2+\frac{1}{2}N(2-x^2)+\frac{1}{2}M(6x-x^2)+x\r]y''(x)
+\l[1-(M+1)x\r]y'(x)+{}\nn
& &+\frac{1}{120}n\l[MN(n^2-1)(n+2)(2n+1)+10Nn(n^2-1)+60M(n+1)+120\r]y(x)=0,\n\eea
for the polynomials
$$y(x):=L_n^{0,M,N}(x)=A_0L_n(x)+A_1L_n'(x)+A_2L_n''(x)$$
with
$$L_n(x):=L_n^{(0)}(x)=\sum_{k=0}^n\frac{(-1)^k}{k!}{n \ch k}x^k,\;\ndots$$
and
$$\l\{\ba{l}\ds A_0=1+Mn+\frac{1}{6}Nn(n-1)(2n-1)+\frac{1}{12}MNn^2(n^2-1)\\ \\
\ds A_1=M+Nn(n-1)+\frac{1}{3}MNn(n^2-1)\\ \\
\ds A_2=Nn+\frac{1}{2}MNn(n+1).\ea\r.$$

The results obtained in these three special cases give rise to the following
conjecture.

\vspace{5mm}

{\bf Conjecture.} {\sl If $\a$ is a nonnegative integer, the polynomials
$\set{\gL}$ satisfy a differential equation of formal order $4\a+10$ which
is of the form
\bea & &M\sum_{i=0}^{2\a+4}\a_i(x)y^{(i)}(x)+
N\sum_{i=0}^{2\a+8}\b_i(x)y^{(i)}(x)+{}\nn
& &\hspace{1cm}{}+MN\sum_{i=0}^{4\a+10}\c_i(x)y^{(i)}(x)+
xy''(x)+(\a+1-x)y'(x)+ny(x)=0,\n\eea
where the coefficients $\l\{\a_i(x)\r\}_{i=1}^{2\a+4}$,
$\l\{\b_i(x)\r\}_{i=1}^{2\a+8}$ and $\l\{\c_i(x)\r\}_{i=1}^{4\a+10}$ are
polynomials independent of $n$ which satisfy
$$\sum_{i=1}^{2\a+4}\a_i(x)=\sum_{i=1}^{2\a+8}\b_i(x)=
\sum_{i=1}^{4\a+10}\c_i(x)=0.$$}

\section{The generalized Jacobi polynomials}

In this section we will deal with the problem of finding a differential
equation for the generalized Jacobi polynomials $\set{\GP}$.

Since the well-known second order differential equation for the classical
Jacobi polynomials $\set{\P}$ is given by
$$(1-x^2)y''(x)+\l[\b-\a-(\a+\b+2)x\r]y'(x)+n(n+\a+\b+1)y(x)=0,$$
it is clear that we look for a differential equation of the form
\bea & &M\sum_{i=0}^{\infty}a_i(x)y^{(i)}(x)+N\sum_{i=0}^{\infty}b_i(x)y^{(i)}(x)
+MN\sum_{i=0}^{\infty}c_i(x)y^{(i)}(x)+{}\nn
& &\hspace{1cm}{}+(1-x^2)y''(x)
+\l[\b-\a-(\a+\b+2)x\r]y'(x)+n(n+\a+\b+1)y(x)=0,\n\eea
where
$$y(x)=P_n^{\a,\b,M,N}(x)=A_0P_n^{(\a,\b)}(x)+\l[A_1(1-x)-A_2(1+x)\r]\frac{d}{dx}P_n^{(\a,\b)}(x)$$
and
$$\l\{\ba{l}\ds A_0=1+M\frac{\ds {n+\b \ch n-1}{n+\a+\b+1 \ch n}}{\ds {n+\a \ch n}}
+N\frac{\ds {n+\a \ch n-1}{n+\a+\b+1 \ch n}}{\ds {n+\b \ch n}}+{}\\  \\
\ds\hspace{7cm}{}+MN\frac{(\a+\b+2)^2}{(\a+1)(\b+1)}{n+\a+\b+1 \ch n-1}^2\\ \\
\ds A_1=\frac{M}{(\a+\b+1)}\frac{\ds {n+\b \ch n}{n+\a+\b \ch n}}{\ds {n+\a \ch n}}
+\frac{MN}{(\a+1)}{n+\a+\b \ch n-1}{n+\a+\b+1 \ch n}\\ \\
\ds A_2=\frac{N}{(\a+\b+1)}\frac{\ds {n+\a \ch n}{n+\a+\b \ch n}}{\ds {n+\b \ch n}}
+\frac{MN}{(\b+1)}{n+\a+\b \ch n-1}{n+\a+\b+1 \ch n}.\ea\r.$$
Here we used the same definition as in \cite{Koorn}, but in a slightly different
notation.

In this case the differential equation will be unique in its general form if
it exists. We introduce the notation
$$\l\{\ba{l}a_0(x)=a_0(n,\a,\b,x)\\ \\
a_i(x)=a_i(\a,\b,x),\;i=1,2,3,\ldots,\ea\r.$$
$$\l\{\ba{l}b_0(x)=b_0(n,\a,\b,x)\\ \\
b_i(x)=b_i(\a,\b,x),\;i=1,2,3,\ldots\ea\r.$$
and
$$\l\{\ba{l}c_0(x)=c_0(n,\a,\b,x)\\ \\
c_i(x)=c_i(\a,\b,x),\;i=1,2,3,\ldots.\ea\r.$$
Since the polynomials $\set{\GP}$ satisfy the symmetry relation
$$P_n^{\a,\b,M,N}(-x)=(-1)^nP_n^{\b,\a,N,M}(x)$$
we have
$$\l\{\ba{l}a_0(n,\a,\b,-x)=b_0(n,\b,\a,x)\\  \\
a_i(\a,\b,-x)=(-1)^ib_i(\b,\a,x),\;i=1,2,3,\ldots\ea\r.$$
and
$$\l\{\ba{l}c_0(n,\a,\b,-x)=c_0(n,\b,\a,x)\\  \\
c_i(\a,\b,-x)=(-1)^ic_i(\b,\a,x),\;i=1,2,3,\ldots.\ea\r.$$

The general form is still an open problem, but the symmetric case $\b=\a$
and $N=M$ turns out to be much less difficult. If $\b=\a$ and $N=M$ we can
choose $c_i(x)=0$ for all $i=0,1,2,\ldots$.

In fact we have shown that there exists a differential equation of the form
\be\la{dvSGP}M\sum_{i=0}^{\infty}a_i(x)y^{(i)}(x)+
(1-x^2)y''(x)-2(\a+1)xy'(x)+n(n+2\a+1)y(x)=0,\ee
where
\be\la{defSGP}y(x)=\SGP=C_0\SP-C_1x\frac{d}{dx}\SP\ee
and
\be\la{CC}\l\{\ba{l}\ds C_0=1+M\frac{2n}{(\a+1)}{n+2\a+1 \ch n}
+4M^2{n+2\a+1 \ch n-1}^2\\
\\
\ds C_1=\frac{2M}{(2\a+1)}{n+2\a \ch n}
+\frac{2M^2}{(\a+1)}{n+2\a \ch n-1}{n+2\a+1 \ch n}.\ea\r.\ee

This differential equation turns out not to be unique. We write
$$\l\{\ba{l}a_0(x):=a_0(n,\a,x),\;\ndots\\  \\
a_i(x):=a_i(\a,x),\;i=1,2,3,\ldots.\ea\r.$$

If we substitute (\ref{defSGP}) and (\ref{CC}) in the differential equation
(\ref{dvSGP}) then we finally find three equations for the coefficients
$\l\{a_i(x)\r\}_{i=0}^{\infty}$ which are equivalent to the following two~:
$$\sum_{i=0}^{\infty}a_i(x)D^i\SP=\frac{4}{(2\a+1)}{n+2\a \ch n}
\frac{d^2}{dx^2}\SP$$
and
$$\sum_{i=0}^{\infty}ia_i(x)D^i\SP+x\sum_{i=0}^{\infty}a_i(x)D^{i+1}\SP
=4{n+2\a+1 \ch n-1}\frac{d^2}{dx^2}\SP.$$

We introduce the notation
$$\l\{\ba{l}a_0(n,\a,x)=a_0(1,\a,x)b_0(n,\a,x)+c_0(n,\a,x),\;\ndots\\ \\
a_i(\a,x)=a_0(1,\a,x)b_i(\a,x)+c_i(\a,x),\;i=1,2,3,\ldots.\ea\r.$$

Now $a_0(1,\a,x)$ is arbitrary, but $\l\{b_i(x)\r\}_{i=0}^{\infty}$ and
$\l\{c_i(x)\r\}_{i=0}^{\infty}$ are uniquely determined. The explicit form
of these coefficients is given in the two theorems below.

\st{The polynomials $\set{\SGP}$ satisfy the linear infinite order differential
equation given by~:
$$\sum_{i=0}^{\infty}b_i(x)y^{(i)}(x)=0,$$
where
$$\l\{\ba{l}\ds b_0(x):=b_0(n,\a,x)=\frac{1}{2}\l[1-(-1)^n\r],\;\ndots\\
\\
\ds b_i(x):=b_i(\a,x)=\frac{2^{i-1}}{i!}(-x)^i,\;i=1,2,3,\ldots.\ea\r.$$}

Now we come to our main result~:

\st{The polynomials $\set{\SGP}$ satisfy the differential equation given by~:
$$M\sum_{i=0}^{\infty}c_i(x)y^{(i)}(x)+(1-x^2)y''(x)-2(\a+1)xy'(x)+n(n+2\a+1)y(x)=0$$
where the coefficients $\l\{c_i(x)\r\}_{i=0}^{\infty}$ are defined by
$$c_0(x):=c_0(n,\a,x)=4(2\a+3){n+2\a+2 \ch n-2},\;\ndots$$
and
$$c_i(x)=(2\a+3)(1-x^2)c_i^*(x),\;i=1,2,3,\ldots,$$
where
$$\l\{\ba{l}\ds c_1^*(x):=c_1^*(\a,x)=0\\ \\
\ds c_i^*(x):=c_i^*(\a,x)=\frac{2^i}{i!}\sum_{k=0}^{i-2}
{\a+1 \ch i-k-2}{i-2\a-5 \ch k}\l(\frac{1-x}{2}\r)^k,\\
\hfill i=2,3,4,\ldots.\ea\r.$$}

We remark that this theorem implies that if $\a$ is a nonnegative integer,
the polynomials $\set{\SGP}$ satisfy a linear differential equation of formal
order $2\a+4$.

Finally, we note that since
$$\hyp{2}{1}{-i+1,\a+\frac{5}{2}-i}{\frac{1}{2}}{1}=\frac{(-\a-2+i)_{i-1}}
{(\frac{1}{2})_{i-1}},\;i=1,2,3,\ldots$$
and
$$\hyp{2}{1}{-i+1,\a+\frac{5}{2}-i}{\frac{3}{2}}{1}=\frac{(-\a-1+i)_{i-1}}
{(\frac{3}{2})_{i-1}},\;i=1,2,3,\ldots$$
we have
\bea \sum_{i=1}^{\infty}c_{2i}^*(\a,-1)&=&
\sum_{i=1}^{\infty}\frac{4}{(\frac{1}{2})_i(1)_i4^i}\frac{(-\a-1)_{i-1}}
{(i-1)!}\frac{(-\a-2+i)_{i-1}}{(\frac{1}{2})_{i-1}}\nn
&=&2\sum_{i=0}^{\infty}\frac{(-\a-1)_i(-\a-1+i)_i}{i!(\frac{3}{2})_i(2)_i
(\frac{1}{2})_i}(\frac{1}{4})^i
=2\sum_{i=0}^{\infty}\frac{(-\a-1)_{2i}}{(2i)!(3)_{2i}}2^{2i}\n\eea
and
\bea \sum_{i=0}^{\infty}c_{2i+1}^*(\a,-1)&=&
-\sum_{i=1}^{\infty}\frac{8(\a+1)}{(1)_i(\frac{3}{2})_i4^i}\frac{(-\a)_{i-1}}
{(i-1)!}\frac{(-\a-1+i)_{i-1}}{(\frac{3}{2})_{i-1}}\nn
&=&-\frac{4}{3}(\a+1)\sum_{i=0}^{\infty}\frac{(-\a)_i(-\a+i)_i}{i!(2)_i
(\frac{5}{2})_i(\frac{3}{2})_i}(\frac{1}{4})^i\nn
&=&-\frac{4}{3}(\a+1)\sum_{i=0}^{\infty}\frac{(-\a)_{2i}}{(2i+1)!(4)_{2i}}2^{2i}
=2\sum_{i=0}^{\infty}\frac{(-\a-1)_{2i+1}}{(2i+1)!(3)_{2i+1}}2^{2i+1}.\n\eea
In the same way we find
$$\sum_{i=1}^{\infty}c_{2i}^*(\a,1)=
2\sum_{i=0}^{\infty}\frac{(-\a-1)_{2i}}{(2i)!(3)_{2i}}2^{2i}$$
and
$$\sum_{i=0}^{\infty}c_{2i+1}^*(\a,1)=
-2\sum_{i=0}^{\infty}\frac{(-\a-1)_{2i+1}}{(2i+1)!(3)_{2i+1}}2^{2i+1}.$$
This implies that
\bea \sum_{i=1}^{\infty}c_i^*(\a,-1)&=&\sum_{i=1}^{\infty}c_{2i}^*(\a,-1)+
\sum_{i=0}^{\infty}c_{2i+1}^*(\a,-1)\nn
&=&2\sum_{i=0}^{\infty}\frac{(-\a-1)_i}{i!(3)_i}2^i=2\;\hyp{1}{1}{-\a-1}{3}{2}\n\eea
and
\bea \sum_{i=1}^{\infty}c_i^*(\a,1)&=&\sum_{i=1}^{\infty}c_{2i}^*(\a,1)+
\sum_{i=0}^{\infty}c_{2i+1}^*(\a,1)\nn
&=&2\sum_{i=0}^{\infty}\frac{(-\a-1)_i}{i!(3)_i}(-2)^i=2\;\hyp{1}{1}{-\a-1}{3}{-2}.\n\eea

\vspace{5mm}

{\bf Remark.} The proof of theorem 5 will be given in
a forthcoming paper \cite{SGP}. Further we remark that this result proves
one of the conjectures by L.L.~Little\-john and W.N.~Everitt given in
\cite{Conj}.

\vspace{5mm}

{\bf Acknowledgement.} The author wishes to thank Professors Desmond Evans,
Norrie Everitt and Lance Littlejohn for their invitation to come to Cardiff,
where the formulas of theorem 5 were found.

Further he wishes to thank his uncle Jan Koekoek for his valuable remarks
after reading the first version of this paper.

\end{document}